%% file: main.tex
\documentclass[11pt,a4paper]{article}

\usepackage{amsmath}
\usepackage{pdflscape}
\usepackage{amssymb}
\usepackage{amsthm}
\usepackage{bm}
\usepackage{authblk}
\usepackage{cite}
\usepackage[english]{babel}
\usepackage[T1]{fontenc}
\usepackage{graphicx}
\usepackage{mathtools}
\usepackage{empheq}
\usepackage{subfig}
\usepackage[utf8]{inputenc}
\usepackage{pgf}
\usepackage{tikz}
\usepackage{url}
\usepackage{bm}
\usepackage{arydshln}
\usepackage{multirow}
\usepackage[draft]{todonotes}   
\usepackage[font=scriptsize]{caption}
\usepackage{mwe}

\graphicspath{ {./Figures/} }

\newtheorem{lemm}{Lemma}

\theoremstyle{definition}
 
\theoremstyle{definition}

\numberwithin{equation}{section}

\newcommand{\TOL}{\texttt{tol}}

\newcommand{\Ad}{\mathrm{Ad}}
\newcommand{\ad}{\mathrm{ad}}
\newcommand{\morespacearray}{\renewcommand{\arraystretch}{1.4}}

\definecolor{rednode}{rgb}{0.8,0.2,0.3}
\definecolor{greennode}{rgb}{0.2,0.8,0.3}
\definecolor{bluenode}{rgb}{0.3,0.2,0.8}
\title{Variable stepsize commutator free Lie group integrators}
\author{Charles Curry \and Brynjulf Owren}
\begin{document}
\maketitle

\begin{abstract}
We introduce variable stepsize commutator free Lie group integrators, where the error control is achieved using embedded Runge-Kutta pairs.
\end{abstract}

\input{introduction}
\input{cfree}
\input{practical}
\input{Vanderpol}
\input{HeavyTop}
\bibliographystyle{plain}
\bibliography{cfree_adaptive}

\end{document}

%% file: introduction.tex
\section{Introduction}
Off the shelf computer software for the numerical solution of ordinary differential equations usually comes with built-in error control and variable step size. Typically, such codes compute an estimate of the local truncation error in each step. This estimate is then compared to a user specified tolerance, \TOL. If its norm is smaller than \TOL, the step will be accepted and a step size to be used in the succeeding step is computed. If the estimate is larger than \TOL, the step will be rejected and attempted again with a new reduced step size. Several adjustments to this procedure can be made in order to improve the behaviour of the integrator. For instance one can seek to avoid an excessive number of rejected steps, or an oscillating behaviour of the sequence of step sizes.

For Runge-Kutta schemes, the most popular method for obtaining an error estimate is by using the method of embedded pairs, see for instance the monograph \cite[p. 164--172]{hairer93sod1} for details.
For a given initial point $y_n$  and a proposed step size $h_n$,  $y_{n+1}$ is computed with a method of convergence order $p$, together with an auxiliary second approximation $\hat{y}_{n+1}$ of convergence order $\hat{p}\neq p$. The resulting pair is usually labelled $p(\hat{p})$ in the literature.
Historically, it was customary to use $\hat{p}>p$ and take $y_{n+1}-\hat{y}_{n+1}$ as an estimate of the error in the lower order approximation $y_{n+1}$, a popular scheme of this form is the Runge--Kutta-Fehlberg 4(5) method. Yet, it seems unfortunate to not make use of the more accurate approximation, and it looks like the preferable choice these days is  to use $\hat{p}<p$ but adjust the step size changing formula in such a way that the global error behaves proportionally to the user given tolerance $\TOL$. This is done in several codes, such as the Dormand-Prince 5(4) method \cite{dormand80afo} implemented in the MATLAB solver \texttt{ode45} \cite{shampine97tmo}.
An important feature of embedded Runge--Kutta pairs is that they share the same 
 internal stages. In this way, the cost of computing with the pair of method is not much higher than computing only the principal approximation $y_{n+1}$ in each step.
\newline
In the last few decades a subfield of numerical analysis called geometric integration has been established, see \cite{hairer10gni} for an exhaustive account. The main purpose is to develop numerical methods that preserve certain underlying geometric structures of the differential equation. Examples of such structures are symplecticity, volume, reversibility, and first integrals. In this paper we consider schemes which were designed for the situation where the differential equation has a natural formulation by means of a Lie group action on a smooth manifold, these schemes are called Lie group integrators, \cite{celledoni14ait}, and the particular type considered here are called commutator-free Lie group integrators \cite{celledoni03cfl}.
Many of the aforementioned methods have primarily been implemented with constant step size.
In fact, for certain geometric properties of a scheme, such as symplecticity, it is not straightforward to vary the step size in such a way that the property is preserved.
\newline
Turning  now to Lie group integrators, it is for our purpose useful to divide them into two categories. One  consists of schemes which can be interpreted as applying a standard Runge--Kutta scheme to a local coordinate representation of the original differential equation, this category includes the Runge--Kutta--Munthe--Kaas methods. For such methods, error control and variable step size can be implemented simply by applying the technique described above to the coordinate representation. The second category, however, does not have such a natural coordinate representation, the schemes are typically constructed by composing flows of simple vector fields. In this case, embedded pairs can only be derived by solving order conditions for Lie group integrators. For a general and modern exposition of order conditions of such integrators, see for example \cite{lundervold15oas}. A detailed account of the order conditions for commutator-free Lie group integrators was given in \cite{owren06ocf}.
The problem of deriving embedded pairs of efficient commutator-free schemes requires careful consideration. Not only is it important that the stages of the pair can be shared, in order for such schemes to be competitive to other Lie group integrators, they also need to reuse flow calculations to as large an extent as possible.

The aim of this paper is to show how one can derive embedded pairs of commutator-free Runge--Kutta methods with the smallest additional computational cost measured in terms of stages and flow calculations (exponentials) per step. Of course, it by no means clear that this way of measuring efficiency yields the smallest possible global error for a given computational cost. But we believe that our approach is a good starting point that illustrates in a clear way some of the challenges involved in constructing embedded pairs of commutator-free methods.
The rest of this paper is organised as follows: in section~2 we begin by briefly describing the class of schemes we consider, and in particular we review the order conditions which were derived in \cite{owren06ocf}. We then give particular examples of embedded pairs of type 3(2) and 4(3). These schemes are constructed such that the pairs share common stages, and employ the well-known ``First same as last'' (FSAL) property known from classical Runge-Kutta methods. We also design the pair to reuse flow calculations to a largest possible extent.
In the third section, we apply the CF43 method to three examples; first the Euler free rigid body, a standard test case for Lie group integrators. 
As a second example, we apply the method to the stiff Van der Pol equation, where we formulate the problem via the standard matrix-vector action by $GL(2)$ on $\mathbb{R}^2$. Finally, we consider the heavy top, formulated via the coadjoint action of the Lie group SE(3) on the dual $\mathfrak{se}(3)^*$ of its  Lie algebra.

%% file: cfree.tex
\section{Commutator free methods}
We consider the numerical solution of autonomous initial value problems on manifolds. For a given manifold $M$ and vector field $F$, that is a section of the tangent bundle $TM$, a smooth curve $y(t)$ is an integral curve of $F$ if it satisfies
\[
\dot{y} = F(y(t)), \quad t\in(a,b)
\]
For any given starting location $y(a)=p$, the integral curve through $p$ exists for sufficiently small time intervals. In the special case that $M$ is a finite-dimensional vector space (which is identified with $\mathbb{R}^n$), integral curves may be approximated numerically using Runge-Kutta methods, i.e.
\begin{eqnarray*}
k_i &=& f\left(y_n + h\sum_{j=1}^s a^j_i k_j\right),\quad i=1,\ldots,s \\
y_{n+1} &=& y_n + h \sum_{i=1}^s b^i k_i
\end{eqnarray*}
It is not so apparent how to extend such methods to manifolds not possessing a linear structure. We begin by observing that the Euler scheme may be interpreted as a method which first `freezes' the vector field $f$ at the point $p=y_n$ so that $f_p(y)=f_p(y_n)$ for all $y$, and then flows a distance $h$ along the integral curve of the frozen vector field $f_p$, i.e.
\[
y_{n+1} = \exp(hf_p)y_n,
\] 
where $\exp(tf_p)$ is the one-parameter group of diffeomorphisms of $\mathbb{R}^n$ induced by the vector field $f_p$. However, our notion of freezing the vector field at $p$ relies on the identification of the tangent spaces $T_p(M)$ and $T_y(M)$ for any $y$. In the general setting, this requires a notion of parallel transport, which in turn requires a choice of linear connection. In general, if we freeze a vector field $f$ by parallel transporting $f_p(p)\in T_p(M)$ onto $T_y(M)$, the result is not independent of the path connecting $p$ and $y$ along which the transport is performed. For Lie group integrators, the connection is typically chosen to be flat so we can ignore this issue, but more generally assuming $y$ is sufficiently close to $p$, the path could be taken to be the geodesic connecting $y$ and $p$. In any case, the consequence is that $\exp(tf_p).p$ will coincide with the geodesic $\gamma_t$ in the direction of $f(p)$ from $p$.

Two strategies emerge from the above comments when generalizing Runge-Kutta methods to manifolds. The first \cite{crouch93nio,owren99rkm,celledoni03cfl} consists of considering the $k_i$ as vector fields frozen at different points, parallel transporting these to a base tangent space $T_p(M)$, taking a linear combination of the resulting tangent vectors and flowing along the geodesic in this direction. This strategy, when implemented using a Lie group action on a homogeneous space, results in the RKMK family of integrators, and allows the use of standard Runge-Kutta tableaux.

In some situations, it is advantageous to avoid using the linear structure on $T_p(M)$ by moving entirely by composition of exponentials of frozen vector fields. This leads to the commutator-free Runge-Kutta methods:
\begin{eqnarray*}
g_r &=& \exp(h \sum_k \alpha^k_{r,J} f_{g_k})\cdots \exp(h\sum_k\alpha^k_{r,1} f_{g_k}) p,\quad r=1,\ldots,s \\
y_{n+1} &=& \exp(h\sum_k \beta^k_J \,f_{g_k})\cdots\exp(h\sum_k\beta^k_1 \,f_{g_k}) p
\end{eqnarray*}

\subsection{Lie group integrators}
The above methods are in principle very general, but require the computation of geodesics and possibly parallel transport, which is typically impractical. The observation underlying Lie group integrators is that the geodesics $\exp(hf_p)$ on a Lie group or reductive homogeneous space equipped with the canonical connection can be computed using the Lie group exponential. In practice, this typically means a matrix exponential, which can be approximated to machine accuracy with tolerable computational effort. Indeed, we suppose that the ODE we wish to solve may be written in the form
\begin{equation}
y' = \big( \lambda_* f(y)\big)(y),\quad y(0)=p,
\label{eq:mkform}
\end{equation}
where $f:M\rightarrow\mathfrak{g}$, and $\lambda_*:\mathfrak{g}\rightarrow\mathcal{X}(M)$ is the infinitesimal action arising from a group action $\Lambda:G\times M\rightarrow M$ as
\[
\lambda_*(u)(p) = \left.\frac{d}{dt}\right|_{t=0} \Lambda(\exp(tu),p)
\]
The vector field $\lambda_* f$ is then frozen straightforwardly at a point $p$ by taking
\[
\lambda_*f_p(y) = \left.\frac{d}{dt}\right|_{t=0} \Lambda(\exp(tf(p)),y)
\]
The exponentials of these frozen vector fields are seen to obey
\[
\exp(t\lambda_* f_p)y = \Lambda(\exp(tf(p)),y)
\]
As a consequence, commutator-free RK Lie-group integrators take the form
\begin{eqnarray*}
g_r &=& \exp(h \sum_k \alpha^k_{r,J} f_k)\cdots \exp(h\sum_k\alpha^k_{r,1} f_k) \\
f_r &=& f(\Lambda(g_r,p)) \\
y_{n+1} &=& \Lambda\big(\exp(h\sum_k \beta^k_J \,f_k)\cdots\exp(h\sum_k\beta^k_1 \,f_k),p\big)
\end{eqnarray*}

\subsection{Order conditions}

Let $\hat{\varphi}_h:M\rightarrow M$ be the mapping corresponding to taking a single timestep of a given commutator-free RK method, and $\varphi_h$ the mapping which flows for time $h$ along the solution curve through the initial point. The method is said to be of order $p$ if, for all $C^{\infty}$-functions $\psi$ on $M$, we have
\[
\psi(\varphi_h(p))-\psi(\hat{\varphi}_h(p)) = \mathcal{O}(h^{p+1})
\]

A commutator-free RK scheme is in general specified by the collection of coefficients $\alpha^k_{r,j}$,$\beta^k_j$. An order theory comprising systems of algebraic equations in the coefficients to be satisfied to attain a given order was derived in \cite{owren06ocf}, analogous to the theory of order conditions of standard RK methods. We will summarize the results without proof here. First, define the coefficients
\[
a^k_r = \sum_j \alpha^k_{r,j},\quad b^k = \sum_j \beta^k_j,\quad c_k = \sum_i a^i_k
\]
\begin{lemm}
A necessary condition for a commutator-free RK scheme with coefficients $\alpha^k_{r,j}$,$\beta^k_j$ to have order $p$ is that the associated $a^k_r$, $b^k$ are the coefficients of a standard Runge-Kutta method of order $p$ or greater.
\end{lemm}

To attain a method of order $1$ or $2$, the above conditions are also sufficient. In particular, we can take $J=1$, such that there is only one exponential computed at each stage. On the other hand, we must satisfy an additional condition to attain order 3:

\begin{lemm}
There are no order 3 methods employing exclusively $J=1$. Given that $J=2$, a sufficient condition to attain order 3 is that the coefficients $a$ and $b$ form a classical RK3 method, and in addition
\[
\sum_k \beta_1^k c_k + \frac{1}{2} \beta_2^k = \frac{1}{3},
\] 
\end{lemm}

We see that we only require two exponentials for computation of $y_{n+1}$; one exponential suffices for the computation of $g_r$ and $f_r$, i.e. we can take $\alpha^k_{r,1}=a^k_r$. The situation changes for methods of order 4:

\begin{lemm}
Given $J=2$, the order 4 conditions are the classical conditions together with the non-classical order 3 condition and the following:
\begin{eqnarray*}
\sum_k \beta_1^k c_k + \frac{1}{3} \beta_2^k &=& \frac{1}{4} \\
\sum_k \beta_1^k c_k^2 + \frac{1}{3} \beta_2^k &=& \frac{1}{6} \\
\sum_{j,k} \beta_1^k a_k^j c_j + \frac{1}{6} \sum_k \beta_2^k &=& \frac{1}{12} \\
\sum_{i,j} b^i c_i \alpha^i_{j,1} c_j +  \sum_{i,j,k} b^i a_i^j c_i \alpha^i_{j,2} c_j &=& \frac{1}{12}
\end{eqnarray*}
\end{lemm}
To satisfy these conditions, it satisfies to take two exponentials in one of the four intermediate stages, two in the final stage, and one exponential in other stages. We illustrate this by giving a sample tableau, displaying the coefficients of a method given in \cite{owren06ocf} which extends the classical RK4 method:
\[
\morespacearray
\begin{array}{r|rrrrr}
0 &&&& \\
\frac{1}{2} & \frac{1}{2} &&& \\
\frac{1}{2} & 0 & \frac{1}{2} && \\
\multirow{2}{*}{1} & \frac{1}{2} & 0 & 0 & \\
& -\frac{1}{2} & 0 & 1 & \\
\hline
& \frac{1}{4} & \frac{1}{6} & \frac{1}{6} & -\frac{1}{12} \\
& -\frac{1}{12} & \frac{1}{6} & \frac{1}{6} & \frac{1}{4}
\end{array}
\]
This is a four-stage method, using two exponentials for the fourth stage, where the tableau displayed above is such that $\alpha_{4,1}=(\frac{1}{2},0,0)$, and $\alpha_{4,2}=(-\frac{1}{2},0,1)$.

To obtain a method of order 5, at least three exponentials are required for the final stage. There results a large nonlinear system of algebraic equations which has so far proved resistant to all attempts to attain a solution. Indeed, the construction of commutator-free methods of order 5 or higher remains an open problem.

\subsection{Embedded pairs reusing exponentials}

Classical Runge-Kutta schemes are typically implemented as an embedded pair, i.e. two different sets of coefficients $b^i, \hat{b}^i$ are given such that the associated approximations $y_n,\hat{y}_n$ are of different orders, typically $\hat{p}=p-1$. In this section, we extend this idea to commutator free RK methods, devising schemes with coefficients $\beta^i_j,\hat{\beta}^i_j$ to create embedded pairs of orders differing by 1.

In practice, it is important to note that the dominant computational cost in the implementation of a commutator free method is likely to be the evaluation of Lie group (matrix) exponentials. Each {\em{unique}} horizontal row in the tableau generally requires the evaluation of one exponential. The observation of Owren \cite{owren06ocf} was that it is possible to design tableaux where certain rows coincide, such as the commutator free RK4 method above where the second and fourth rows are identical. This consideration looms large in the implementation of embedded pairs, where there is even greater opportunity for reuse of exponentials. A general CF3 scheme requires at least 3 exponentials, whilst a CF4 scheme requires 5. Our main achievement is the construction of CF32 and CF43 schemes using only one additional exponential and function evaluation compared to the constant stepsize scheme of the same order. 

\subsubsection{CF32}
The order conditions to obtain a commutator free method of second order coincide with the classical RK2 conditions. As in the classical case the construction of a RK32 pair is always possible from an RK3 scheme, there is no difficulty in constructing a commutator free RK32 pair. In general, this will require one extra exponential and function evaluation. Suppose the tableau is constructed as follows
\[
\morespacearray
\begin{array}{r|ccc}
0 &&& \\
c_2 & a_{21} && \\
c_3 & a_{31} & a_{32} & \\
\hline
\multirow{2}{*}{$y$} & \beta^1_1 & \beta^2_1 & \beta^3_1  \\
& \beta^1_2 & \beta^2_2 & \beta^3_2  \\
\hline
\hat{y} & a_{31} & a_{32} & 0,
\end{array}
\quad\mathrm{or}
\quad
\begin{array}{r|cccc}
0 &&&& \\
c_2 & a_{21} &&& \\
c_3 & a_{31} & a_{32} && \\
\multirow{2}{*}{$1$} & \beta^1_1 & \beta^2_1 & \beta^3_1 &  \\
& \beta^1_2 & \beta^2_2 & \beta^3_2  & \\
\hline
\multirow{2}{*}{$y$} & \beta^1_1 & \beta^2_1 & \beta^3_1 &  \\
& \beta^1_2 & \beta^2_2 & \beta^3_2  & \\
\hline
\hat{y} & a_{31} & a_{32} & 0 & 0,
\end{array}
\]
implying a reuse of the third stage and hence saving an exponential. The two schemes above are identical, but the right hand tableau is of a more general type that allows $\hat{y}$ to depend on the value of $f(y)$. Schemes of this type have the FSAL (first same as last) property, as the first function evaluation at the next step coincides with the function evaluation $f(y)$ assuming the step is accepted. We will henceforth write all tableaux in this form, and abbreviate the row for $y$ as FSAL.

The order 2 conditions become $a_{31}+a_{32}=1$ and $a_{21}a_{32}=\frac{1}{2}$, hence any 3-stage commutator free RK3 scheme obeying
\[
a_{32}=\frac{1}{2a_{21}},\quad a_{31}=\frac{2a_{21}-1}{2a_{21}}
\]
would admit an embedded pair without the need for an extra exponential. Unfortunately, such a scheme does not exist, as is shown readily using symbolic computation software. We therefore focus on constructing commutator free RK3 schemes reusing an exponential, a topic as yet unexplored in the literature. All of the results stated in the remainder of this section were proven using symbolic computation software.

Let $\omega$ be a root of $36z^2 + (9a-30)z + 3a + 1$. The general tableau of an RK32 scheme reusing the third stage in the second row of the fourth is given in Table~\ref{TableFSALgen32}.
\begin{table}
{\footnotesize
\centering 
\[
\morespacearray
\begin{array}{r|cccc}
0 &&&& \\
a & a &&& \\
\frac{6\omega - 1}{3} & \frac{6a\omega -3\omega -a}{3a} & \frac{\omega}{a} && \\
\multirow{2}{*}{$1$} & \frac{-12a\omega-18\omega-11a+11}{6} & \frac{18\omega(1-a) - 11}{6a(3a+1)} & \frac{-12\omega +3a -10}{2(3a+1)}& \\
& \frac{6a\omega -3\omega -a}{3a} & \frac{\omega}{a} &&
\\
\hline
y & \mathrm{FSAL} &&& \\
\hline
\hat{y} &&&&
\end{array}
\]}
{\footnotesize
\[
\morespacearray
\begin{array}{r|rrrr}
0 &&&& \\
\frac{1}{3} & \frac{1}{3} &&& \\
1 & -1 & 2 && \\
\multirow{2}{*}{$1$} & 1 & -\frac{5}{4} & \frac{1}{4} & \\
& -1 & 2 & 0 & \\
\hline
y & \mathrm{FSAL} &&& \\
\hline
\hat{y} & 0 & \frac{3}{4} & 0 & \frac{1}{4}
\end{array}
\qquad\qquad
\begin{array}{r|rrrr}
0 &&&& \\
\frac{1}{3} & \frac{1}{3} &&& \\
-\frac{1}{6} & -\frac{5}{12} & \frac{1}{4} && \\
\multirow{2}{*}{$1$} & -\frac{37}{12} & \frac{9}{4} & 2 & \\
& -\frac{5}{12} & \frac{1}{4} && \\
\hline
y & \mathrm{FSAL} &&& \\
\hline
\hat{y} & 0 & \frac{3}{4} & 0 & \frac{1}{4}
\end{array}
\]
}
\caption{In the first row a we give a general FSAL commutator-free RK-scheme of order 3(2) reusing the exponential of the third stage in the second row of  the fourth stage.  $\omega$ is a root of the polynomial $36z^2 + (9a-30)z + 3a + 1$.
In the second row two concrete examples with rational coefficients are given.}
\label{TableFSALgen32}
\end{table}

Once $a$ is chosen and a consistent choice of root for $\omega(a)$ is made, there then follows a system of two linear equations for the four coefficients of $\hat{y}$, in general resulting in a two-parameter family once the rest of the tableau is set. In general, the values will be irrational, but we highlight that some choices of $a$ give rational coefficients, for instance $a=-\frac{1}{3},\frac{1}{3},\frac{2}{3},\frac{7}{9},\frac{101}{9},\frac{34}{3},\frac{35}{3}$ etc. In general these can be found by setting the discriminant $81a^2 - 972a + 756$ to be a square. We give two samples of rational tableaux from the case $a=\frac{1}{3}$ in the second row of Table~\ref{TableFSALgen32}.
We can also reuse the third stage in the first row of the fourth; indeed this may be preferable for some integrators as it permits the storage of the action of the exponential on $y_0$ rather than the exponential itself. Let $\nu$ be the root of $36z^2 + (9a-6)z -3a + 1$. The general tableau for such a scheme is found in Table~\ref{cfree32r3inst41}.
\begin{table}
\centering
{\footnotesize
\[
\morespacearray
\begin{array}{r|cccc}
0 &&&& \\
a & a &&& \\
\frac{6\nu-1}{3}& \frac{6a\nu-3\nu+a}{3a} & \frac{\nu}{a} && \\[2mm]
\multirow{2}{*}{$1$} & \frac{6a\nu-3\nu+a}{3a} & \frac{\nu}{a} && \\ 
&\frac{-12a\nu-6\nu+a+1}{6a} & \frac{-18a\nu-6\nu + 1}{6a(3a-1)} & \frac{12\nu +3a -2}{2(3a-1)}&
\\
\hline
y & \mathrm{FSAL} &&& \\
\hline
\hat{y} &&&&
\end{array}
\]}
\caption{Commutator-free scheme of order 3(2) reusing the exponential of the third stage in the first row of the fourth stage.}
\label{cfree32r3inst41}
\end{table}
There are also one-parameter families of RK32 schemes reusing the second stage, either in the first or second row of $y$. Fix $a$ and let $\gamma$ be a root of $4a(3a-1)z^2 +4(3a-1)z + 3$ and $\delta$ be a root of $4az^2 + (12a-2)z + 9a+6$. The respective tableaux are given in Table~\ref{cfree32r2inupdate}.
\begin{table}
{\footnotesize
\[
\morespacearray
\begin{array}{r|cccc}
0 &&&& \\
\frac{1}{3} & \frac{1}{3} &&& \\
c_3 & a & (\frac{2}{3}-2a)(2a\gamma + 1) && \\[2mm]
\multirow{2}{*}{$1$} & \frac{1}{3} & 0 & 0 & \\
& (3a-1)\gamma + \frac{2}{3} & -3a\gamma & \gamma & \\
\hline
y & \mathrm{FSAL} &&& \\
\hline
\hat{y} & & & & 
\end{array}
\qquad
\begin{array}{r|cccc}
0 &&&& \\
-\frac{1}{3} & -\frac{1}{3} &&& \\
c_3& -\frac{2a\delta}{3}-2a & a && \\[2mm]
\multirow{2}{*}{$1$} & \frac{-6a\delta +8a+3}{6a} & \delta & -\frac{1}{2a} & \\
& -\frac{1}{3} & 0 &0 & \\
\hline
y & \mathrm{FSAL} &&& \\
\hline
\hat{y} & & &  & 
\end{array}
\]}
\caption{One parameter families of commuttator-free schemes reusing the exponential of the second stage in the first row (left) or the second row (right)}
\label{cfree32r2inupdate}
\end{table}
The discriminants of the equations for $\gamma$ and $\delta$ are $16(1-3a)$ and $4(1-36a)$ respectively, which allows for easy generation of rational versions of the above schemes should this be desired, for instance $a=0$ works in both cases.
\subsubsection{CF43}
In the classical case, there are no RK43 pairs with only 4 stages, so it is essential to use a 5-stage FSAL scheme to attain RK43 at minimal cost. This is no longer true for commutator free methods, as it is possible to take a single CF4 method together with $\beta,\hat{\beta}$ such that the underlying $b=\hat{b}$, but that only $\beta$ fulfills all of the non-classical order four conditions. Nonetheless, we will preserve the greater generality afforded by the FSAL schemes in our presentation.

An order three method requires two exponentials in the final computation, so in the generic case a CF43 pair requires two extra exponentials compared to the CF4 case. It is possible to reduce the extra cost to one by reusing an exponential, but the full tableau must be generated with this in mind, as a generic CF4 scheme does not permit a CF43 pair which reuses an exponential. In general, it is possible to reuse two exponentials (but not more) in a CF43 tableau, one above the lines (in the $\alpha$s) and one below the lines (i.e., generating $\hat{y}$). Once the reuse pattern has been specified, the top of the tableau is fixed, and a one parameter family of $\hat{\beta}$ coefficients is typically admitted. Note that patterns involving reuse of the second stage do not give CF43 pairs. In contrast to the CF32 case, there are no rational CF43 pairs reusing the optimal number of exponentials. We give the exact form of one of the tableaux in Table~\ref{cf43exact}, but typically print only floating point forms; more accurate descriptions of the coefficients are available from the authors on request.

\begin{table}
{\footnotesize
\[
\morespacearray
\begin{array}{r|ccccc}
0 &&&&& \\
c_2 & p_1(\omega) &&&& \\
c_3      &p_2(\omega)
      &p_3(\omega) &&& \\[2mm]
\multirow{2}{*}{$1$}      &p_2(\omega)
    &p_3(\omega) &&& \\
   & p_4(\omega) & p_5(\omega)&p_6(\omega) && \\[2mm]
\multirow{2}{*}{$1$} & p_7(\omega) &
p_8(\omega) &
p_9(\omega) & \frac{\omega}{2}& \\
&-\frac13 p_7(\omega)&
p_{10}(\omega) &
p_{11}(\omega)& \frac{-3\omega}{2} & \\
\hline
y & \mathrm{FSAL} &&&& \\
\hline
\hat{y} 
&p_{12}(\omega) & p_{13}(\omega)&p_{14}(\omega) & 0 & 0 \\
& \mathrm{one\, parameter \,family} & & && \\
\end{array}
\]
\begin{align*}
 p_1(\omega) &= \frac12 (7-288\omega^4-36\omega^3+48\omega^2+17\omega) \\
 p_2(\omega) &= \frac{1}{268}(-389+31824 \omega^4+10962\omega^3-3651 \omega^2-2027 \omega) \\
 p_3(\omega)&=\frac{1}{268}(54-2880 \omega^4-2520 \omega^3+234 \omega^2+553 \omega) \\
 p_4(\omega)&=\frac{-51696\omega^4-13878\omega^3+7557\omega^2+2285\omega+1244}{804}  \\
 p_5(\omega)&=\frac{-521424\omega^4-323586\omega^3+61119\omega^2+61599\omega+10976}{20100}\\
 p_6(\omega)&=\frac{-5328\omega^4+558\omega^3+93\omega^2-122\omega+47}{300} \\
 p_7(\omega)&= \frac{1008\omega^4-1530\omega^3+501\omega^2-16\omega+229}{536} \\
 p_8(\omega)&=\frac{541872\omega^4+76158\omega^3-84207\omega^2-19972\omega-2703}{40200} \\
 p_9(\omega)&=\frac{-2304\omega^4+144\omega^3+174\omega^2+4\omega+21}{150} \\
 p_{10}(\omega)&=\frac{256752\omega^4+67878\omega^3-170787\omega^2-10852\omega+22877}{40200}\\
p_{11}(\omega)&=\frac{-864\omega^4-396\omega^3+684\omega^2+264\omega+11}{150}  \\
p_{12}(\omega)&= \frac{-51696\omega^4-13878\omega^3+7557\omega^2+2285\omega+1244}{804} \\
p_{13}(\omega)&=\frac{-521424\omega^4-323586\omega^3+61119\omega^2+61599\omega+10976}{20100}  \\
p_{14}(\omega)&= \frac{-5328\omega^4+558\omega^3+93\omega^2-122\omega+47}{300}
\end{align*}
}
\caption{A commutator-free pair of order 4(3). Here $\omega$ is the unique real root of $144z^5 +90z^4 -3z^3 -13z^2 -5z-1$. This is the general form of a 4(3)-pair reusing the second stage in the first part of the third stage, and reusing the second part of the third stage in the first part of the $\hat{y}$ computation.}
\label{cf43exact}
\end{table}
The same tableau is given below as a decimal approximation:
{\footnotesize
\[
\morespacearray
\begin{array}{r|rrrrr}
0 &&&&& \\
4.785707347 & 4.785707347 &&&& \\
.8093268944 & .7701000600 &  .03922683443 &&& \\[2mm]
 \multirow{2}{*}{$1$}& .7701000600 & .03922683443 & 0 && \\
& .6195164818 & .06934556872  & -.4981889449   && \\[2mm]
\multirow{2}{*}{$1$} & .4211354919 & -.005776103764 & -.1381183969  & .2227590088 & \\
& -.1403784973 & .006491728470 & 1.302163795  & -.6682770264 & \\
\hline
y & \mathrm{FSAL} &&&& \\
\hline
\hat{y} & $.6195164818$ & $0.06934556872$ & $-.4981889449$ && \\ \hdashline
& -.158427746\,a      &  0.008517658\,a & a & -.618653791\,a  & -.231436127\,a  \\
&     - 0.075415454 & -0.082788288       &    & + .5828295568  & + .3847010797
\end{array}
\]
}

Another possible choice of reuse pattern leads to the next tableau, where the one parameter family for $\beta_2$ has been set such that $\beta^3_2=0$.
{\footnotesize
\[
\begin{array}{r|ccccc}
0 &&&&& \\
 & .67104050 &&&& \\
 & 2.547687640 & -1.355037274 &&& \\
 & 2.547687640 & -1.355037274 & 0 && \\
& -.21944181 & -0.0735967 & .1003880 && \\
1 & .324015249 & .15832891 & -.21057643 & .2282322824 & \\
& -.108005081 & .84426683 & .44843513 & -.6846968472 & \\
\hline
y & \mathrm{FSAL} &&&& \\
\hline
& -.21944181 & -0.0735967 & .1003880 & 0 & 0 \\
& .45603817 & .93310478 & 0 & -.2660264 & 0.06953371 \\
\end{array}
\]
}
Perhaps the optimal CF43 method in terms of reuse pattern is the following, in which the reuse always occurs in the first part of a split stage, and the method is a true 4-stage method not using FSAL. It is the only CF43 scheme with these properties:
{\footnotesize
\[
\begin{array}{r|cccc}
0 &&&& \\
 & 1.351207192 &&& \\
 & 0.5 & 0.097900176 && \\
 & 0.5 & 0.097900176 && \\
 & 7.900943678 & 2.989500877 & -10.48834473 & \\
\hline
y & .301574869 & -0.054881885 & .238291289 & 0.01501572796  \\
& -.1005249562 & .1005249562 & .5450471839 & -0.04504718389  \\
\hline
\hat{y}  & 0.5 & 0.097900176 & 0 & 0  \\
& -.2989500877 & -0.0522571042 & .783338473 & -0.03003145592 \\
\end{array}
\]}

%% file: practical.tex
\section{Practical implementation}
We give an outline of how the above methods are implemented in practice. First we show how an embedded pair allows for automatic step size control, following \cite{hairer93sod1}. Indeed, suppose we have chosen an initial step size $h$, and obtain approximate solutions $y_1$ and $\hat{y}_1$. If $y$ takes values in a normed space, we can consider $y-\hat{y}$ as an estimate of the error, and aim to ensure
\[
||y_1-\hat{y}_1||<sc,\quad sc=Atol + \max(||y_0||,||y_1||)\cdot Rtol,
\]
for some user-specified absolute and relative tolerances $Atol,Rtol$. If $y$ takes values in a manifold equipped with only a metric, the notion of relative error does not make sense, but we can still aim to ensure that $d(y_1,\hat{y}_1)<sc=Atol$. In either case, let $err=d(y_1,\hat{y}_1)/sc$. In general, for an embedded pair of order $p(p-1)$, the optimal step size is $h_{opt}=h\cdot err^{-\frac{1}{p}}$. In practice, it is usual to dampen the fluctuations in $h$, a typically procedure is to let
\[
h_{new} = h\cdot \min(facmax,\max(facmin,fac\cdot err^{-\frac{1}{p}})),
\]
for some $fac<1$ and appropriate $facmin,facmax$. The above discussion is in some sense only rigorous when the lower order approximation $\hat{y}$ is used to continue the integration, as only then is $d(y,\hat{y})$ a good measure of the local error. Nonetheless, practice has shown that it is usually better to continue with the higher order integrator, a procedure known as local extrapolation.

We have not addressed automatic selection of the initial step size. For problems in vector spaces, this is typically performed by constructing estimates of the derivatives of vector field $F$, using some combination of function evaluations and taking small step(s) with the Euler scheme, see \cite{hairer93sod1}. Similar ideas may be employed for Lie group integrators, but require modification as derivatives cannot be approximated so simply using evaluations of vector fields at different points due to the difficulty of identifying nearby tangent spaces.

\subsection{Free rigid body}

We now give a simple illustrative example of Lie group integration. Euler's equations for the body angular momentum of a free rigid body in a reference frame parallel to the principal axes of inertia take the form
\begin{equation}\label{frb}
\frac{d\xi}{dt} = -mI^{-1}\hat{\xi}.\xi,
\end{equation}
where $m$ is the body mass, $I$ is the diagonal inertia tensor, and we use the standard convention of the hat-map
$$
\xi=\left(\begin{array}{c}\xi_1\\ \xi_2\\ \xi_3\end{array}\right)\quad\Rightarrow\quad
\hat{\xi} =\left(\begin{array}{ccc}0 & -\xi_3 & \xi_2\\ \xi_3 & 0 & -\xi_1\\ -\xi_2 & \xi_1 & 0\end{array}\right)
$$
The quantity $|\xi|^2$ is readily shown to be conserved, i.e. equation $\ref{frb}$ describes the evolution of $\xi\in S^2$. The equations allow for the immediate application of Lie group integrators, indeed the right hand side is of the form $\lambda_* f(\xi)\xi$, where the manifold $M=S^2$, $f:S^2\rightarrow \mathfrak{so}(3)$ is the hat map, and the infinitesimal action is that of $\mathfrak{so}(3)$ on $S^2$ by matrix multiplication. The associated group action $\Lambda$ is $SO(3)$ acting by matrix multiplication. In this context, the commutator-free RK Lie-group integrator given by the coefficients $\alpha,\beta$ is therefore
\begin{eqnarray*}
g_r &=& \exp(h \sum_k \alpha^k_{r,J} f_k)\cdots \exp(h\sum_k\alpha^k_{r,1} f_k) \\
f_r &=& \widehat{g_r\cdot \xi_n} \\
\xi_{n+1} &=& \exp(h\sum_k \beta^k_J \,f_k)\cdots\exp(h\sum_k\beta^k_1 \,f_k)\cdot \xi_n,
\end{eqnarray*}
where $g_r\in SO(3)$, $f_r\in\mathfrak{so}(3)$, and $\xi_n\in S^2$.
\begin{figure}
\centering
\includegraphics[width=\textwidth]{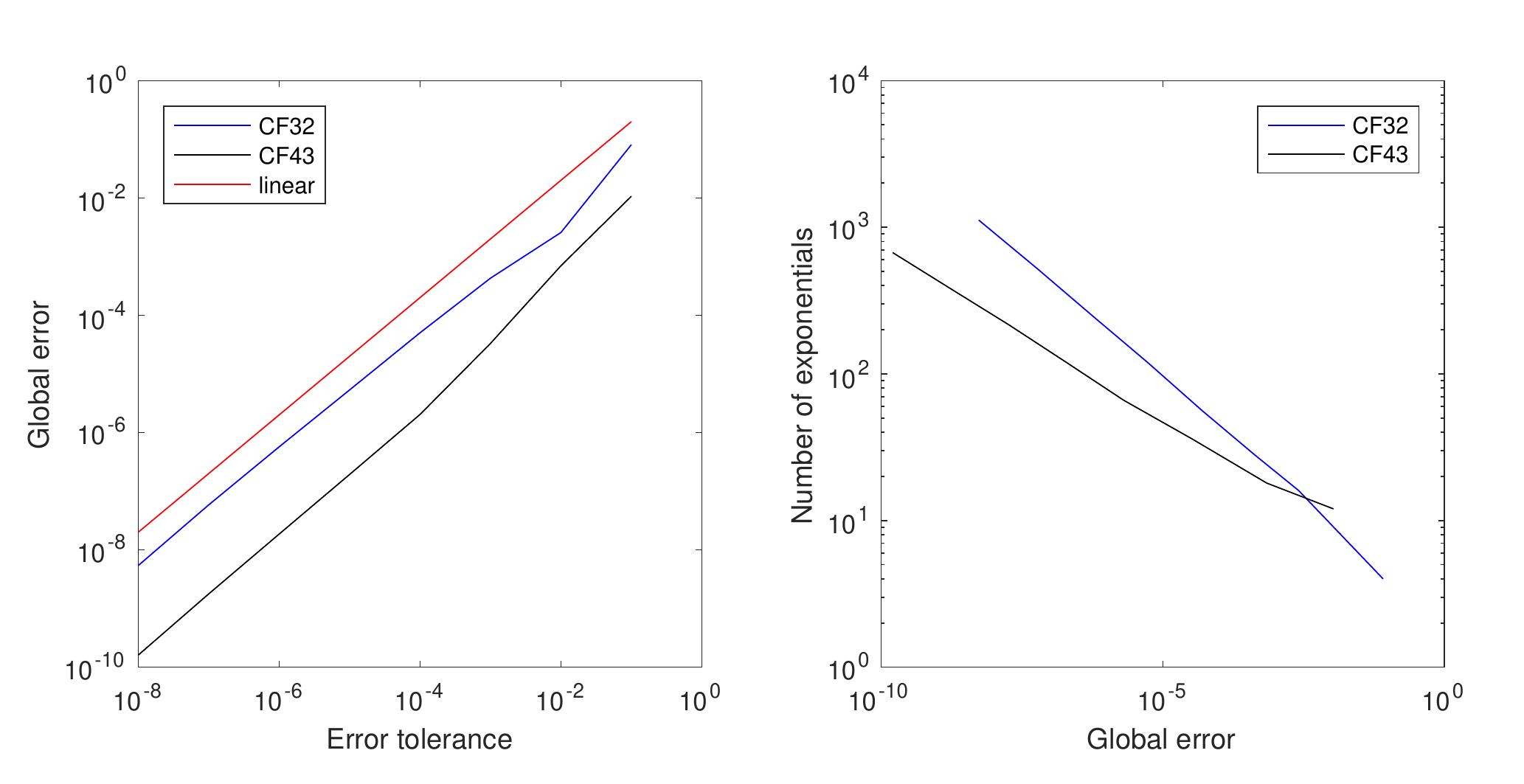}
\caption{The CF32 pair (bottom left, Table~\ref{TableFSALgen32}) and the CF43 pair (Table~\ref{cf43exact}) applied to the free rigid body equations with $I=(1,2,5)^T$ and random initial data on $S^2$ on the time interval $[0,2]$. The graph on the left shows the linear relationship between the prescribed error tolerance and the observed global error, whilst the right graph compares the computational cost (measured in number of exponentials) against global error.
}
\label{figure:vdp1}
\end{figure}

%% file: Vanderpol.tex
\subsection{The Van der Pol oscillator}
We consider the following non-conservative oscillator with non-linear damping, formulated as a scalar second order differential equation
\begin{equation} \label{vanderpol}
\ddot{x} - \mu (1-x^2)\,\dot{x} + x = 0
\end{equation}
where $\mu$ is a parameter that affects the stiffness of the system. One can rephrase this problem in the form \eqref{eq:mkform} using the simple matrix-times-vector action of the Lie group $GL(2)$ on $\mathbb{R}^2\backslash\{0\}$, 
\begin{equation} \label{vanderpol_sys}
\frac{\mathrm{d}}{\mathrm{dt}}\left(\begin{array}{c} x\\ \dot{x}\end{array}\right)
=\left(\begin{array}{cc} 0 & 1 \\ -1 & \mu(1-x^2)\end{array}\right)\,\left(\begin{array}{c} x\\ \dot{x}\end{array}\right)
\end{equation}
The corresponding Lie group integrator can be interpreted as an exponential integrator in the sense defined for instance in
\cite{berland05bsa}. An explicit Lie group integrator cannot be expected to work well for stiff problems, but we believe it is still of interest to observe how the new embedded pair of commutator-free schemes behaves through the ``needle'' of the Van der Pol oscillator.
\begin{figure}
\centering
\includegraphics[width=\textwidth]{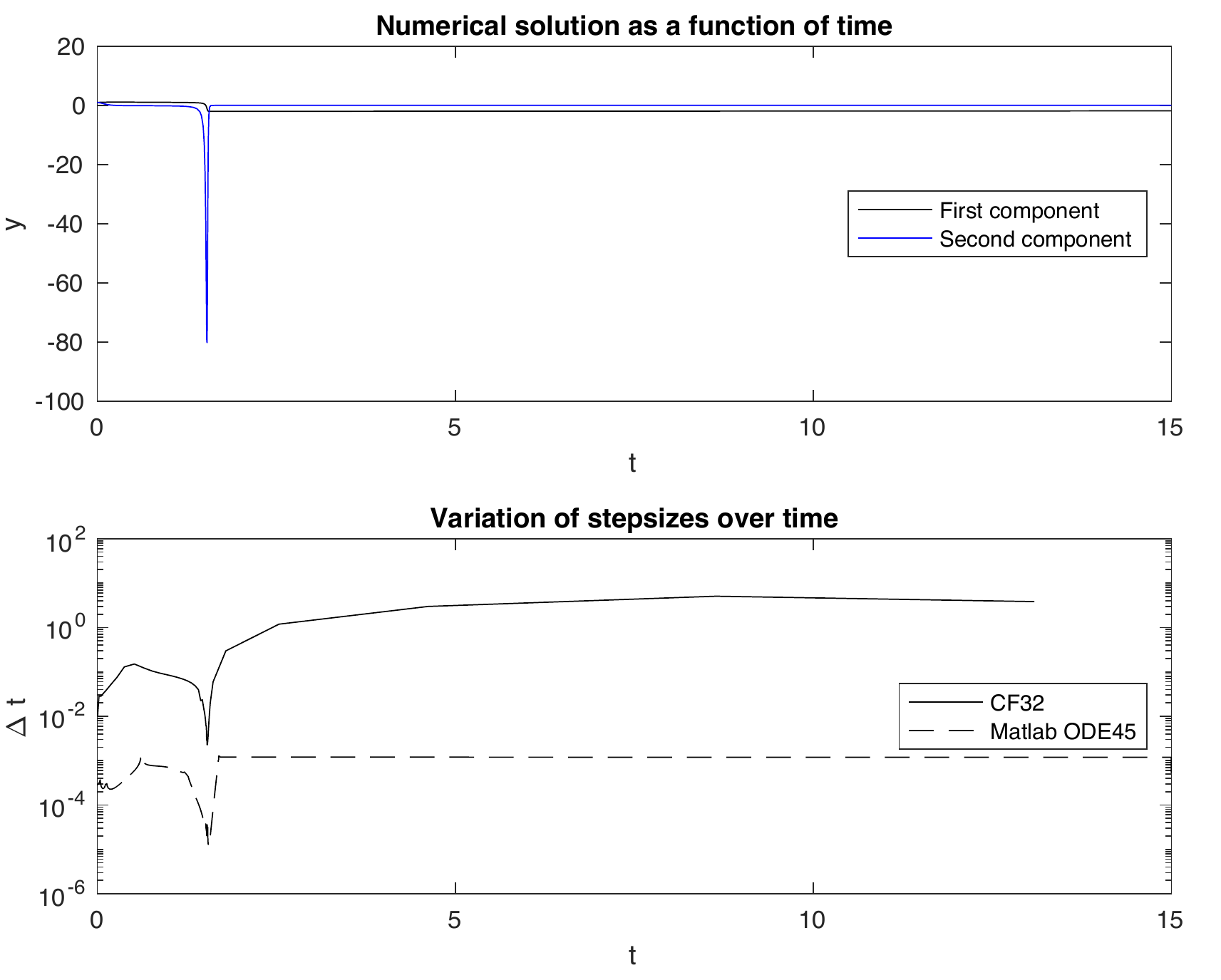}
\caption{The CF32 pair (bottom left, Table~\ref{TableFSALgen32}) applied to the Van der Pol equation with $\mu=60$ and $y_0=(1,1)^T$ on the time interval $[0,15]$. The figure shows the numerical solution (top) and the step size sequences selected by the method (bottom). The dashed line shows for comparison the step size sequence used by the Matlab solver ODE45. Both solvers used a tolerance $\mathtt{TOL}=10^{-3}$.
}
\label{figure:vdp1}
\end{figure}
We have implemented the CF32 scheme given in the bottom left tableau of Table~\ref{TableFSALgen32} and applied it in a variable stepsize fashion to the Van der Pol oscillator \eqref{vanderpol_sys} setting $\mu=60$ and $y_0=(1,1)^T$ in all the experiments. In Figure~\ref{figure:vdp1} we show the two components of the solution computed by CF32 in the top graph, and note in particular the sharp downward spike (``needle")
for the second component in the approximate interval $t\in[1.4,1.56]$. The relative and absolute tolerances were both set to $10^{-3}$ in this experiment.
The stepsizes chosen by the CF32 schemes are shown in the bottom graph of Figure~\ref{figure:vdp1}, and one can see how the stepsizes are reduced through the needle. For comparison, we also show the stepsizes used by the builtin Matlab solver ODE45 which is based on the Dormand-Prince embedded Runge--Kutta pair \cite{dormand80afo}. We observe that the new CF32 solver behaves similarly to the Dormand-Prince scheme except that the former takes larger steps, this might be expected due to the fact that it computes (exact) matrix exponentials. For matrices with eigenvalues whose real parts tend to $-\infty$, such exact exponentials are bounded as opposed to their explicit Runge--Kutta counterparts which use polynomial approximations to the exponential map.
\begin{figure}
\includegraphics[width=\textwidth]{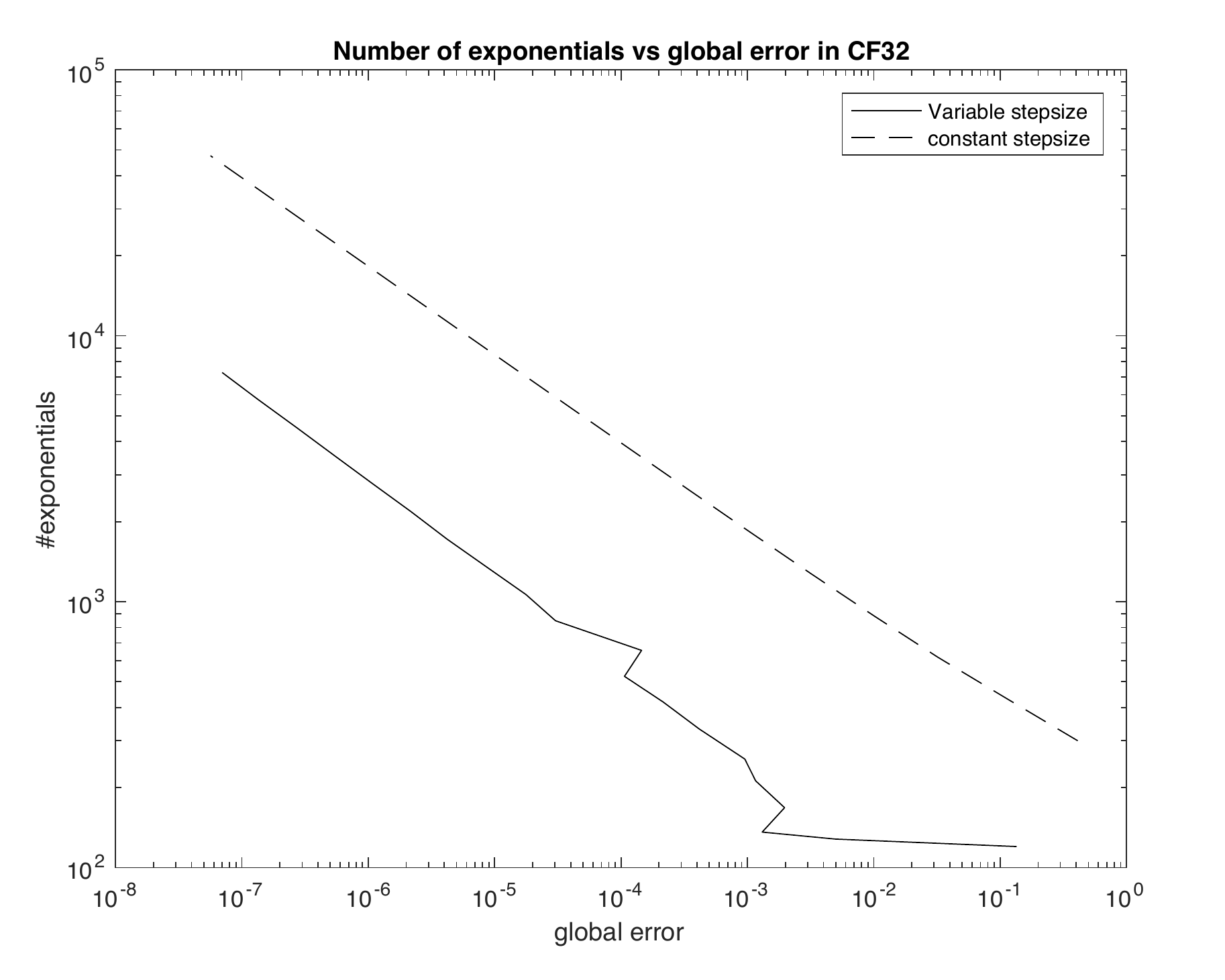}
\caption{Numerical integration of the Van der Pol equation with $\mu=60$, $y_0=(1,1)^T$.
The horizontal axis is the global error at $t=1.6$, and the vertical axis shows the number of exponentials that were computed. The solid line is the variable stepsize method, and the dashed line is the same third order method applied with constant stepsize.
 }
 \label{figure:vdp2}
\end{figure}
In Figure~\ref{figure:vdp2} we visualise the difference between constant and variable stepsize for the problem \eqref{vanderpol_sys}.
By running the variable stepsize code for a number of different choices for tolerances, we compute the global error at a fixed point $t=1.6$, just after passing through the needle. The exact solution has been approximated various different ways, one was to run various builtin Matlab solvers with strict tolerances.  As a measure for the efficiency of the integrator, we have computed  the number of exponential calcualtions used by the code to obtain a prescribed global error. The dashed line in Figure~\ref{figure:vdp2} shows the result for the constant stepsize version of the third order method in our embedded pair, and the solid line shows the result for the variable stepsize solver. The cost ratio between the constant and variable stepsize methods depend on the chosen global error, but for instance to obtain a global error of $10^{-5}$ the constant stepsize integrator needs approximately 6.5 times as many exponential calculations as the variable stepsize method.

%% file: HeavyTop.tex
\subsection{The heavy top}
Mathematical models for the heavy top can be found in many text books, see for instance
\cite{marsden99itm,holm11gmp}. Lie group integrators were applied to this problem in \cite{marthinsen97soo}. The heavy top is a rigid body, but because of the gravitational forces, it is not invariant under the action of $SO(3)$ and the dynamics can therefore not be reduced to $\mathfrak{so}(3)^*$ as the free rigid body.
There is a smaller symmetry group $S^1$ corresponding to rotation about the vertical axis and it turns out that the system can be formulated on $\mathfrak{se}(3)^*$, the dual of the Lie algebra of the special Euclidean group $SE(3)$ consisting of translations and rotations in 3-space. According to \cite{holm11gmp} the equations can be written in the form
\begin{align}
\begin{split}  \label{heavytop}
\dot{\bm{\mu}} &= \bm{\mu}\times \mathbb{I}^{-1}\bm{\mu} + \bm{\beta}\times mg\bm{\chi} \\
\dot{\bm{\beta}} &= \bm{\beta}\times \mathbb{I}^{-1}\bm{\mu}
\end{split}
\end{align}
Here $\bm{\mu}$ is the body angular momentum and $\bm{\beta}$ is the vertical direction as seen from the rotating body, more precisely $\bm{\beta}=R^T\bm{e}_3$ where $R$ is the attitude matrix of the top.
$\bm{\chi}$ is the unit vector in the direction from the fixed point to the center of mass of the heavy top, $m$ is the mass, $g$ the constant of gravity, and $\mathbb{I}$ is the inertia tensor.

In this example, the coadjoint orbits are preserved, and for this reason it is natural to invoke Lie group integrators via the (right) coadjoint action of $SE(3)$ on $\mathfrak{se}(3)^*$. It is convenient to take elements of both $\mathfrak{se}(3)$ and $\mathfrak{se}(3)^*$ to be vectors in 
$\mathbb{R}^3\times\mathbb{R}^3$. Similarly, elements of $SE(3)$ are represented as pairs $(g,\bm{u})$ where $g$ could be an orthogonal $3\times 3$-matrix an $\bm{u}\in\mathbb{R}^3$. 

The coadjoint action is the map $\Lambda: SE(3)\times \mathfrak{se}^*(3)\rightarrow \mathfrak{se}(3)^*$
$$
\Lambda((g,\bm{u}),(\bm{\mu},\bm{\beta}))=\Ad^*_{(g,\bm{u})}(\bm{\mu},\bm{\beta})=
 (g^T(\bm{\mu}-\bm{u}\times\bm{\beta}),g^T\bm{\beta})
$$
whose infinitesimal generator is the map $\lambda_*:\mathfrak{se}(3)\rightarrow\mathcal{X}(\mathfrak{se}(3)^*)$
$$
 \lambda_*(\bm{\xi},\bm{u})(\bm{\mu},\bm{\beta})=\ad^*_{(\bm{\xi},\bm{u})}(\bm{\mu},\bm{\beta})
=(-\bm{\xi}\times \bm{\mu}-\bm{u}\times\bm{\beta},-\bm{\xi}\times\bm{\beta})
$$
For the heavy top equations \eqref{heavytop} we now have
$$
(\dot{\bm{\mu}},\dot{\bm{\beta}})=\lambda_*(\mathbb{I}^{-1}\bm{\mu},mg\bm{\chi})(\bm{\mu},\bm{\beta})
$$
so that with reference to \eqref{eq:mkform} we have
$$
      f(\bm{\mu},\bm{\beta})  = (\mathbb{I}^{-1}\bm{\mu},mg\bm{\chi})
$$

The key properties of the group action are summarised in Table~\ref{heavytop-props}.
\begin{table}
\centering
{\footnotesize
\begin{tabular}{|l|l|l|}
\hline 
Lie group   & $SE(3)\cong SO(3)\ltimes\mathbb{R}^3$ \\ 
Group product $SE(3)$ & $(g,\bm{u})\cdot (h,\bm{v}) = (g\cdot h,g\cdot \bm{v} + \bm{u})$  \\ 
Inverse $SE(3)$ &   $(g,\bm{u})^{-1} =(g^{-1},-g^{-1}\bm{u})$ \\ \hline
Lie algebra & $\mathfrak{se}(3)\cong \mathfrak{so}(3)\ltimes\mathbb{R}^3$ \\ 
Lie bracket $\mathfrak{se}(3)$ & $[(\bm{\xi},\bm{u}),(\bm{\eta},\bm{v})]=(\bm{\xi}\times\bm{\eta}, \bm{\xi}\times \bm{v}-\bm{\eta}\times \bm{u})$ \\ \hline
Manifold & $\mathfrak{se(3)^*} \cong \mathbb{R}^3 \times \mathbb{R}^3$ \\ \hline
\begin{minipage}{3.7cm} Coadjoint action by $SE(3)$ on $\mathfrak{se}(3)^*$\end{minipage}  & $\begin{array}{lcl}(g,\bm{u})\cdot(\bm{\mu},\bm{\beta}) &=& \Ad_{(g,\bm{u})}^*(\bm{\mu},\bm{\beta})\\
&=&(g^T(\bm{\mu}-\bm{u}\times\bm{\beta}),g^T\bm{\beta})\end{array}$ \\ \hline
\begin{minipage}{3.8cm} Infinitesimal generator of the action\end{minipage} & $\begin{array}{lcl}\lambda_*(\bm{\xi},\bm{u})(\bm{\mu},\bm{\beta})&=&\ad^*_{(\bm{\xi},\bm{u})}(\bm{\mu},\bm{\beta})\\
&=&(-\bm{\xi}\times \bm{\mu}-\bm{u}\times\bm{\beta},-\bm{\xi}\times\bm{\beta})\end{array}$ \\ \hline
Exponential map & $\exp(t(\bm{\xi},\bm{u}))=(\exp(t\hat{\bm{\xi}}), \frac{\exp(t\hat{\bm{\xi}})-I}{t\hat{\bm{\xi}}}\cdot t\bm{u})$ \\ \hline
\end{tabular}
}
\caption{The main properties of the group action for the heavy top equations}
\label{heavytop-props}
\end{table}

\begin{figure}
\centering
\includegraphics[width=\textwidth]{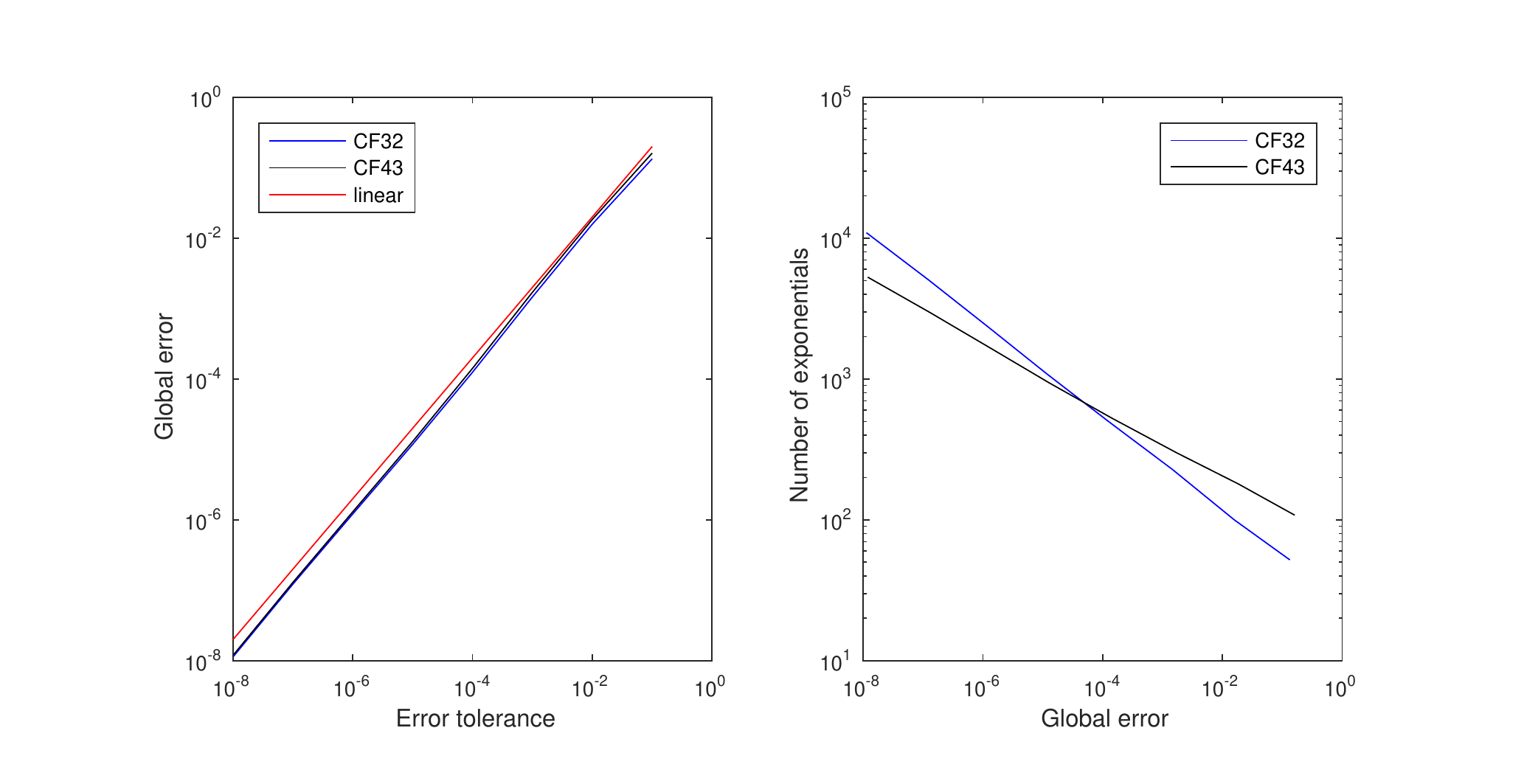}
\caption{The CF32 pair (bottom left, Table~\ref{TableFSALgen32}) and the CF43 pair (Table~\ref{cf43exact}) applied to the Kovalevskaya top on the time interval $[0,2]$. The graph on the left shows the linear relationship between the prescribed error tolerance and the observed global error, whilst the right graph compares the computational cost (measured in number of exponentials) against global error.
}
\label{figure:vdp1}
\end{figure}

%% file: main.bbl
\begin{thebibliography}{10}

\bibitem{berland05bsa}
H.~Berland, B.~Owren, and B.~Skaflestad.
\newblock {$B$}-series and order conditions for exponential integrators.
\newblock {\em SIAM J. Numer. Anal.}, 43(4):1715--1727, 2005.

\bibitem{celledoni03cfl}
E~Celledoni, A~Marthinsen, and B~Owren.
\newblock Commutator-free {L}ie group methods.
\newblock {\em FUTURE GENERATION COMPUTER SYSTEMS}, 19(3):341--352, APR 2003.

\bibitem{celledoni14ait}
E.~Celledoni, H.~Marthinsen, and B.~Owren.
\newblock An introduction to {L}ie group integrators---basics, new developments
  and applications.
\newblock {\em J. Comput. Phys.}, 257(part B):1040--1061, 2014.

\bibitem{crouch93nio}
P.~E. Crouch and R.~Grossman.
\newblock Numerical integration of ordinary differential equations on
  manifolds.
\newblock {\em J. Nonlinear Sci.}, 3(1):1--33, 1993.

\bibitem{dormand80afo}
J.~R. Dormand and P.~J. Prince.
\newblock A family of embedded {R}unge-{K}utta formulae.
\newblock {\em J. Comput. Appl. Math.}, 6(1):19--26, 1980.

\bibitem{hairer10gni}
E.~Hairer, Ch. Lubich, and G.~Wanner.
\newblock {\em Geometric numerical integration}, volume~31 of {\em Springer
  Series in Computational Mathematics}.
\newblock Springer, Heidelberg, 2010.
\newblock Structure-preserving algorithms for ordinary differential equations,
  Reprint of the second (2006) edition.

\bibitem{hairer93sod1}
E.~Hairer, S.~P. N{\o}rsett, and G.~Wanner.
\newblock {\em Solving ordinary differential equations. {I}}, volume~8 of {\em
  Springer Series in Computational Mathematics}.
\newblock Springer-Verlag, Berlin, second edition, 1993.
\newblock Nonstiff problems.

\bibitem{holm11gmp}
D.~D. Holm.
\newblock {\em Geometric mechanics. {P}art {II}. {R}otating, translating and
  rolling}.
\newblock Imperial College Press, London, second edition, 2011.

\bibitem{lundervold15oas}
A.~Lundervold and H.~Z. Munthe-Kaas.
\newblock On algebraic structures of numerical integration on vector spaces and
  manifolds.
\newblock In {\em Fa{\`a} di {B}runo {H}opf algebras, {D}yson-{S}chwinger
  equations, and {L}ie-{B}utcher series}, volume~21 of {\em IRMA Lect. Math.
  Theor. Phys.}, pages 219--263. Eur. Math. Soc., Z{\"u}rich, 2015.

\bibitem{marsden99itm}
J.~E. Marsden and T.~S. Ratiu.
\newblock {\em Introduction to mechanics and symmetry}, volume~17 of {\em Texts
  in Applied Mathematics}.
\newblock Springer-Verlag, New York, second edition, 1999.
\newblock A basic exposition of classical mechanical systems.

\bibitem{marthinsen97soo}
A.~Marthinsen, H.~Z. Munthe-Kaas, and B.~Owren.
\newblock Simulation of ordinary differential equations on manifolds: some
  numerical experiments and verifications.
\newblock {\em Modeling, Identification and Control}, 18(1):75--88, 1997.

\bibitem{owren99rkm}
B.~Owren and A.~Marthinsen.
\newblock Runge-{K}utta methods adapted to manifolds and based on rigid frames.
\newblock {\em BIT}, 39(1):116--142, 1999.

\bibitem{owren06ocf}
Brynjulf Owren.
\newblock Order conditions for commutator-free {L}ie group methods.
\newblock {\em J. Phys. A}, 39(19):5585--5599, 2006.

\bibitem{shampine97tmo}
L.~F. Shampine and M.~W. Reichelt.
\newblock The {MATLAB} {ODE} suite.
\newblock {\em SIAM J. Sci. Comput.}, 18(1):1--22, 1997.
\newblock Dedicated to C. William Gear on the occasion of his 60th birthday.

\end{thebibliography}
